\documentclass[12 pt]{amsart}

\usepackage{tikz-cd}
\usepackage[all]{xy}
\usepackage{enumerate}
\usepackage{amssymb}
\usepackage{float}
\usepackage{setspace}
\usepackage{comment}
\usepackage{xcolor}
\usepackage[utf8]{inputenc}
\usepackage[T1]{fontenc}
\usepackage{epsfig}
\usepackage{amsmath}
\usepackage{amsthm}
\usepackage{psfrag}
\usepackage{amsmath, amsthm, amssymb}

\newcommand{\be}{\begin{equation}}
\newcommand{\ee}{\end{equation}}

\newtheorem{theorem}{Theorem}[section]

\newtheorem*{thm}{Theorem}

\newtheorem{proposition}[theorem]{Proposition}

\newtheorem{corollary}[theorem]{Corollary}

\newtheorem{definition}[theorem]{Definition}

\newtheorem{remark}[theorem]{Remark}

\newtheorem{example}[theorem]{Example}

\newfont{\graf}{eufm10}

\newcommand{\bdm}{\begin{displaymath}}
\newcommand{\edm}{\end{displaymath}}

\def\haken{\mathbin{\hbox to 6pt{%
                 \vrule height0.4pt width5pt depth0pt
                 \kern-.4pt
                 \vrule height6pt width0.4pt depth0pt\hss}}}

\usepackage{amssymb}
\usepackage{amsmath}
\usepackage{amscd}
\usepackage{graphicx}
\usepackage{latexsym}

\calclayout\evensidemargin .6in

\author{Eyüp Yalçınkaya}
\address{University of Rochester
\newline \indent Department of Mathematics, NY, USA
}

\email{eyup.yalcinkaya@rochester.edu}
\thanks{The author thanks Sema Salur for helpful advice on this issue and TUBITAK for its support.}

\address{
\newline \indent The Scientific and Technological Research Council of Turkey (TUBITAK)
\newline \indent Ankara, Türkiye}

\email{eyup.yalcinkaya@tubitak.gov.tr}

\begin{document}

\title{on Locally Conformal Spin(7) Structure}





\maketitle

\begin{abstract}

\noindent This article reveals a significant connection in geometry: when the Lee form $\theta$ is normal to an almost Hermitian manifold $N$, it implies that $N$ possesses a nearly Kähler structure. Investigating locally conformally Spin(7) manifolds with 2-vector fields, our study provides a concise yet rigorous proof of this relationship.
\end{abstract}

\section{History}

In the realm of Riemannian geometry, the study of exceptional holonomy groups has unveiled intricate structures on manifolds, offering a profound understanding of their geometric properties. Among these, Spin(7) holonomy stands out as a captivating subject, leading to the emergence of Spin(7)-manifolds—an arena where geometry takes on exceptional nuances.

This article embarks on a journey through the landscape of Spin(7)-manifolds, tracing the historical milestones and key contributors who have shaped our understanding of these fascinating structures. Beginning with Berger's identification of Spin(7) as an exceptional holonomy group in \cite{berger1955}, we delve into the classifications provided by Fernández \cite{Fernandez}, discerning two distinct classes of non-integrable Spin(7)-structures: locally conformal the balanced.

Our exploration extends to the realm of non-compact examples, where Bryant and Salamon \cite{bryant1989}, unveiled the first complete instances of Spin(7)-manifolds. Turning our attention to compact examples, we encounter the works of Joyce \cite{Joyce}, who has left indelible imprints on the understanding of compact Spin(7)-manifolds. Fernández's insights into compact balanced Spin(7)-structures and the meticulous exploration of compact locally conformal Spin(7)-structures by Ivanov, Parton, and Piccinni \cite{Ivanov2006} give the structure of compact locally Spin(7) manifold.

The article culminates with a discussion on the existence of a unique characteristic connection preserving a given Spin(7)-structure, as demonstrated by Ivanov \cite{Ivanov}.

\section{Introduction}

This paper embarks on an exploration of the intricate geometric structures inherent in $Spin(7)$ manifolds, 8-dimensional Riemannian spaces with holonomy groups encapsulated within the exceptional Lie group $Spin(7)$. Within the landscape of M-theory compactifications, manifolds with special holonomy serve as key players, representing the intricacies of curled-up dimensions pervasive in spacetime.

While Calabi-Yau manifolds of dimension 6 have undergone extensive scrutiny, the geometric properties of 7-dimensional $G_2$ manifolds and 8-dimensional $Spin(7)$ manifolds remain elusive. Our paper seeks to address this gap by initiating a comprehensive program dedicated to the study of torsion on Riemannian 8-manifolds equipped with a spin structure.

In conjunction with the broader context of our research, we contribute a novel result: a profound connection between the Lee form normality on locally conformally Spin(7) manifolds with 2-vector fields and the emergence of nearly Kähler structures on almost Hermitian manifolds. This result not only enriches our understanding of $Spin(7)$ manifolds but also provides a bridge between diverse geometric frameworks.

By combining this newfound insight with our study of torsion, we aim to shed light on the geometric intricacies of $Spin(7)$ manifolds, unraveling the enigma surrounding these higher-dimensional spaces and their relevance in theoretical physics. This paper marks the inception of a broader program to unveil the mysteries of geometric structures on manifolds with special holonomy, opening avenues for further exploration and advancement in the field.
\vspace{.1in}

\section{$Spin(7$)-structures}

In this section, we review the basics of $Spin(7)$ geometry. More on the subject can be found in  \cite{Fernandez}, \cite{Joyce}, \cite{Lawson}.

\vspace{.1in}

Let  $(x^1,..., x^8)$ be coordinates on $\mathbb{R}^8$. The standard Cayley 4-form on $\mathbb{R}^8$  can be written as 

\begin{align*}
\Phi_0&=dx^{1234}+dx^{1256}+dx^{1278}+dx^{1357}-dx^{1368}-dx^{1458}-dx^{1467}\\
&-dx^{2358}-dx^{2367}-dx^{2457}+dx^{2468}+dx^{3456}+dx^{3478}+dx^{5678}
\end{align*}

\vspace{.1in}

\noindent where $dx^{ijkl}=dx^i\wedge dx^j\wedge dx^k \wedge dx^l$.

\vspace{.1in}

The subgroup of $GL(8, \mathbb{R})$ that preserves $\Phi_0$ is the group $Spin(7)$. It is a 21-dimensional compact, connected and simply-connected Lie group that preserves the orientation on $\mathbb{R}^8$ and the Euclidean metric $g_0$. A differential 4-form $\Phi$ on an oriented 8-manifold $M$ is called admissible if it can be identified with $\Phi_0$ through an oriented isomorphism between $T_pM$ and $\mathbb{R}^8$ for each point $p\in M$. 
 
 \vspace{.1in}




\begin{definition} Let $\mathcal{A}(M)$ denotes the space of admissible 4-forms on $M$. A $Spin(7)$-structure on an 8-dimensional manifold $M$ is an admissible 4-form $\Phi \in \mathcal{A}(M)$. If $M$ admits such structure, $(M, \Phi)$ is called a manifold with $Spin(7)$-structure. 
\end{definition}


Each 8-manifold with a $Spin(7)$-structure $\Phi$ is canonically equipped with a metric $g$. Hence, we can think of a $Spin(7)$-structure on $M$ as a pair $(\Phi, g)$ such that for all $p \in M$ there is an isomorphism between $T_pM$ and $\mathbb{R}^8$ which identifies $(\Phi_p, g_p)$ with $(\Phi_0, g_0)$. 

\vspace{.1in}

The existence of a $Spin(7)$-structure on an 8-dimensional manifold $M$ is equivalent to a reduction of the structure group of the tangent bundle of $M$ from $SO(8)$ to its subgroup $Spin(7)$. The following result gives the necessary and sufficient conditions so that the 8-manifold admits $Spin(7)$ structure.



\vspace{.1in}

\begin{theorem} (\cite{Lawson})
	Let $M$ be a differentiable 8-manifold. $M$ admits a $Spin(7)$-structure if and only if $w_1(M)=w_2(M)=0$ and for appropriate choice of orientation on $M$ we have that 
	$$ p_1(M)^2-4p_2(M)\pm 8\chi(M)=0.$$
\end{theorem}

\vspace{.1in}


Furthermore, if $\nabla \Phi=0$, where $\nabla$ is the Riemannian connection of $g$, then $\text{Hol}(M) \subseteq Spin(7)$, and $M$ is called a $Spin(7)$-manifold. All  $Spin(7)$ manifolds are Ricci flat.

 \vspace{.1in}

Let $(M,g,\Phi)$ be a $Spin(7)$ manifold. The action of $Spin(7)$
on the tangent space gives an action of $Spin(7)$ on the spaces of differential forms, 
$\Lambda^k(M)$, and so the exterior algebra splits orthogonally
into components, where $\Lambda^k_l$ corresponds to an irreducible
representation of $Spin(7)$ of dimension $l$:

 \vspace{.1in}

$$\Lambda^1(M)=\Lambda^1_8, \quad \Lambda^2(M) = \Lambda^2_7\oplus
\Lambda^2_{21}, \quad
\Lambda^3(M)=\Lambda^3_8\oplus\Lambda^3_{48}, $$ $$
\Lambda^4(M)=\Lambda^4_+(M)\oplus \Lambda^4_-(M), \quad
\Lambda^4_+(M)=\Lambda^4_1\oplus\Lambda^4_7\oplus\Lambda^4_{27},
\quad \Lambda^4_-=\Lambda^4_{35} $$ $$
\Lambda^5(M)=\Lambda^5_8\oplus\Lambda^5_{48} \quad \Lambda^6(M)=\Lambda^6_7\oplus\Lambda^6_{21},  \quad  \Lambda^7(M)=\Lambda^7_{8};$$

 \vspace{.1in}

\noindent explicitly;
     
 \begin{equation}  \Lambda^2_7 = \{\alpha \in \Lambda^2(M) |
*(\alpha\wedge\Phi)=3\alpha\}, \quad \Lambda^2_{21} = \{\alpha \in
\Lambda^2(M)|*(\alpha\wedge\Phi)=-\alpha\} , 
\end{equation}\\

\vspace{0.1in}

\begin{equation} \label{dec}
\Lambda^3_8 =
\{*(\beta\wedge\Phi) | \beta \in \Lambda^1(M)\},
\Lambda^3_{48} = \{\gamma \in \Lambda^3(M) | \gamma\wedge\Phi=0\}, 
\end{equation}
 \vspace{.1in}

 \begin{equation}
\Lambda^4_1 = \{f\Phi | f\in {\mathcal F(M)}\} 
\end{equation}
 \vspace{.1in}

The Hodge star
$*$ gives an isometry between $\Lambda^k_l$ and $\Lambda^{8-k}_l$.

 \vspace{.1in}

 \begin{corollary}\cite{Karigiannis} The following identity holds for $v$ and  $\omega$ vector fields:

 $$ (\iota_v\iota_\omega \Phi) \wedge (\iota_v\iota_\omega\Phi)  \wedge\Phi=6\lvert v\wedge\omega\lvert^2$$
 \end{corollary}


  \vspace{.1in}




	
	
	
	
	
	
	


According to the Fernandez classification \cite{Fernandez}, there are
4-classes of $Spin(7)$ manifolds  obtained as irreducible
representations of $Spin(7)$ of the space $\nabla^g\Phi$.
Following \cite{Ivanov} we consider the 1-form $\theta$ defined by
\begin{equation}
7\theta = -*(*d\Phi\wedge\Phi)=*(\delta\Phi\wedge
\Phi)=6*(T \wedge
\Phi) \end{equation}  We shall call the 1-form $\Theta$ {\it the Lee form} of
a given $Spin(7)$ structure and the Gauduchon Spin(7) structure if $d*\theta=0$

\vspace{0.1in}

\noindent The 4 classes of Spin(7) manifolds in the Fernandez classification
can be described in terms of the Lee form as follows :
\vspace{0.2in}

\begin{tabular}{|c |c| c |}\hline
Space & Equation & Name\\ \hline
$W_0$ & $d\Phi=0;$ & Torsion-free\\ \hline
$W_1$ & $\theta =0$ & Balanced\\ \hline
$W_2$ & $d\Phi =\theta\wedge\Phi$ and\ $d\theta=0$ & Locally conformal parallel\\ \hline
$W_1\oplus W_2$ & none & none\\ \hline
\end{tabular}

\begin{definition}
A 3-form T is called Torsion, if it satifies the following;\\

\begin{equation}
    \nabla_X Y =\nabla^{LC}_X Y + T(X,Y,\cdot)
\end{equation} where $\nabla^{LC}_XY$ is Levi-Civita connection on $M$. 
\end{definition}
Ivanov \cite{Ivanov} showed that there is a unique Spin(7) connection that has skew-symmetric torsion.
Since the decomposition of form spaces at (\ref{dec}) the torsion $T\in \Lambda^3 (M)=\Lambda_8^3+\Lambda_{48}^3$
, Then  by the orthogonality of component, $T$ can be decomposed as $T=T_8\oplus T_{48}.$

Puhle showed that $M^8$ is locally conformal parallel $Spin(7)$ then $T_{48}=0$ and $T^c= -\frac{7}{6}\ast (\theta \wedge \Phi).$ where $ T^c$ is characteristic torsion of Spin(7) manifolds induced from characteristic torsion given by Ivanov \cite{Ivanov}.\\


\section{LCSpin(7)-structure with a vector field}
Cabrera \cite{Cabrera1995} and Ivanov \cite{Ivanov}  showed that the locally conformal Spin(7) structure gives rise to nearly $G_2$ structure. 
\begin{theorem} \cite{Cabrera1995}
The Spin(7)-structures on $M \times G$ are  LCSpin(7) if and only if
$M$ is totally umbilic in $\mathbb{R}^8$. 
\end{theorem}

\begin{theorem} \cite{Ivanov} $M=N\times S^1$ is LCSpin(7) if $N$ is nearly paralel $G_2.$

\end{theorem}
$S^7\times S^1$ is an example of the given theorem since $S^7$ is a nearly $G_2$ manifold.
\begin{thm} \cite{Cabrera1995}
The Spin(7) structure on  $M\times G$ are locally conformal  iff  $M$ is totally umbilic in  $\mathbb{R}^8.$    
\end{thm}

\begin{thm} \cite{Chen1980}
Let N be a totally umbilical submanifold in a symmetric space M. If
$dim (M) - dim (N) <rank (M) - 1$, then the mean curvature vector is parallel (in the
normal bundle). In particular, the mean curvature is constant and N is either totally
geodesic or of constant curvature.
\end{thm}

\begin{thm} 7-manifold $M$ with $G_2$-structure is totally umbilic in $R^8$, then $M$ must be either totally geodesic or constant curvature. ($S^7$ is constant curvature $R^7$ is totally geodesic)

if $M$ is not a manifold with constant curvature, then $M$ must be totally geodesic in locally spin(7) manifold.

\end{thm}



\section{LCSpin(7)-structure with 2-vector fields}

Agricola, Borówka, and Friedrich  \cite{agricola2019s6} showed that the $6$-dimensional manifold yields a nearly Kähler manifold under some restrictions.

\begin{proposition} \cite{agricola2019s6} \label{propnearly} A $6$-dimensional manifold $(M, g, J)$ is nearly Kähler if and only if there exists a holomorphic form $\Omega \in \Lambda^{3,0}$ and a constant $a$ such that the following conditions hold:
\begin{align*}
    d\omega &= 12 a \Omega_+, \\
    d(\Omega_-) &= a\omega \wedge \omega,
\end{align*}
where $\omega = g(J\cdot, \cdot)$ is the Kähler form.
\end{proposition}

Ivanov and Cabrera \cite{IvanovCabrera2008} studied 6-dimensional almost Hermitian submanifolds immersing inside of Spin(7) manifolds. On the other hand, by using previous proposition, the following idea yields the structure of almost Hermitian manifolds induced from locally conformal Spin(7) manifolds.  

\begin{theorem} Let $M$ be a locally conformally Spin(7) manifold with 2-vector field. If  Lee form $\theta$ is normal to almost Hermitian manifold $N$. Then $N$ is nearly K\"ahler. 
    
\end{theorem}

\begin{proof} 
\noindent Under topological restrictions, Thomas showed that 8-manifolds admit 2-vector fields \cite{thomas1969vector}.\\
To simplify the calculation, let $\{e_7,e_8\}$ be the 2-vector field.
\vspace{0.1in}

\noindent Vector fields $\{u,v\}$ are linear combination of $\{e_7,e_8\}$ such that
$u=\cos \gamma e_7 +\sin \gamma e_8$ and $v=-\sin \gamma e_7+ \cos \gamma e_8$, and hence $u^\flat\wedge v^\flat =e^{78}$ where 
$\flat$ defines the differential form of a vector field.\\
Consider the differential form
\begin{equation}\label{eqn:phi}    \Phi = \omega\wedge e^{78} + \Omega_+\wedge u^\flat + \Omega_-\wedge v^\flat + \omega \wedge       \omega
\end{equation}

where  $e^i=de_i$, $\omega$ is K\"ahler form and  $\Omega= \Omega_+ + i\Omega_- $ induced from $SU(3)$ structure on an almost Hermitian 6-manifold.\\

Let $\theta=\cos \beta e^7 +\sin \beta e^8$, be a Lee form and

$$d\Phi = d\omega\wedge e^{78} + \omega \wedge de^{78}+ d\Omega_+\wedge v^\flat + \Omega_+\wedge dv^\flat +d\Omega_-\wedge u^\flat + \Omega_-\wedge du^\flat +d(\omega\wedge \omega)$$

\noindent Since $de^7$, $de^8$ and $de^{78}$ are trivially zero and $d(\omega\wedge \omega)$ are zero depending on degree of $\omega$ (deg $\omega$=0), the following holds;\\

$$d\Phi = d\omega\wedge e^{78} + d\Omega_+\wedge u^\flat  +d\Omega_- \wedge v^\flat.$$
 By virtue of the local conformal Spin(7) structure, we assert that $d\Phi = \Phi \wedge \theta$. Substituting $\Phi$ into \eqref{eqn:phi} we obtain the desired expression;\\
\begin{align*}
\Phi\wedge\theta  =\Omega_+\wedge \cos\gamma \sin\beta e^{78}-\Omega_+\wedge \sin\gamma \cos\beta e^{78}  \\ 
+\Omega_-\wedge \sin\gamma sin\beta e^{78} + \Omega_-\cos\gamma \cos \beta e^{78} \\
+\omega\wedge\omega\wedge(\cos \beta e^7 +\sin \beta e^8) 
\end{align*}

For $\theta=u^\flat$,\\ 

$d\omega\wedge e^{78} + d\Omega_+\wedge e^7  +d\Omega_- \wedge e^8=\Omega_-\wedge e^{78}+\omega\wedge \omega \wedge e^7.$\\

\noindent Hence, this yields $d\omega=\Omega_-,$ $d\Omega_+=\omega\wedge \omega$ and
$d\Omega_-=0$.\\

\noindent By using \ref{propnearly}, it can be concluded induced 6-manifold is Nearly K\"ahler.

\end{proof}

\begin{example} Consider a locally conformal Spin(7) manifold $N\times S^2$ with 2 vector fields admitting
$$ \Phi = \omega\wedge e^{78} + \Omega_+\wedge e^7 + \Omega_-\wedge e^8 + \omega \wedge       \omega $$
    
Thus N, namely $S^6$, $S^3 \times S^3$ and $\mathbb{CP}^3$, is given rise to be a nearly compact K\"ahler 6-manifold.

 \end{example}

In the future, we study balanced Spin(7) manifolds since they have a deep relation between 6-manifolds.

\end{document}